 \documentclass[final,3p,times,twocolumn]{elsarticle}

\usepackage{amssymb}

\usepackage{algorithm}
\usepackage{algorithmic}

\newtheorem{definition}{Definition}

\newtheorem{lemma}{Lemma}
\newtheorem{remark}{Remark}
\newtheorem{theorem}{Theorem}

\journal{Pattern Recognition}

\begin{document}

\begin{frontmatter}



\title{Generalized singular value thresholding operator to affine matrix rank minimization problem}

 \author[label1]{Angang Cui}
 \ead{cuiangang@163.com}
 \author[label2]{Haiyang Li}
  \ead{fplihaiyang@126.com}
 \author[label3]{Jigen Peng \corref{cor}}
 \cortext[cor]{Corresponding author}
 \ead{jgpengxjtu@126.com}
 \author[label1]{Junxiong Jia}
 \ead{jjx323@xjtu.edu.cn}
 \address[label1]{School of Mathematics and Statistics, Xi'an Jiaotong University, Xi'an, 710049, China.}
 \address[label2]{School of Science, Xi'an Polytechnic University, Xi'an, 710048, China.}
  \address[label3]{School of Mathematics and Information Science, Guangzhou University, Guangzhou, 510006, China.}

\begin{abstract}
It is well known that the affine matrix rank minimization problem is NP-hard and all known algorithms for exactly solving it are doubly exponential in theory
and in practice due to the combinational nature of the rank function. In this paper, a generalized singular value thresholding operator is generated to solve the affine
matrix rank minimization problem. Numerical experiments show that our algorithm performs effectively in finding a low-rank matrix compared with some state-of-art methods.
\end{abstract}

\begin{keyword}
Affine matrix rank minimization problem\sep Generalized thresholding operator\sep Generalized singular value thresholding operator


\MSC 90C26 \sep 90C27 \sep 90C59

\end{keyword}

\end{frontmatter}


\section{Introduction}\label{section1}
The affine matrix rank minimization (AMRM) problem consisting of recovering a low-rank matrix that satisfies a given system of linear equality constraints is an important problem in recent years, and
has attracted much attention in many applications such as machine learning \cite{srebro1}, collaborative filtering in recommender systems \cite{candes2,jan3}, computer vision \cite{hu4}, network localization \cite{ji5}, system identification \cite{fazel6,liu7}, control theory \cite{faze8,faze9}, and so on. A special case of AMRM is the matrix completion (MC) problem \cite{candes10}, it has been applied in the famous Netflix problem \cite{net26} and image inpainting problem \cite{yeg27}. Unfortunately, the problem AMRM is generally NP-hard \cite{recht11} and all known algorithms for exactly solving it are doubly exponential in theory and in practice due to the combinational nature of the rank function.

A popular alternative is the nuclear-norm affine matrix rank minimization (NuAMRM) problem \cite{candes2,fazel6,faze9,candes10,recht11,candes12}.  Recht et al. \cite{recht11} have demonstrated that if a certain restricted isometry property holds for the linear transformation defining the constraints, the minimum rank solution can be recovered by solving the problem NuAMRM. Cai et al.\cite{cai14} considered the  regularization nuclear-norm affine matrix rank minimization (RNuAMRM) problem and a singular value thresholding (SVT) algorithm is proposed to solve this regularization problem. Although there are many theoretical and algorithmic advantages \cite{recht11,cai14,liu13,toh15,ma16} for the convex relaxation problem NuAMRM, it may be suboptimal for recovering a real low-rank matrix and yields a matrix with much higher rank and needs more observations to recover a real low-rank matrix \cite{candes2,cai14}. Moreover, the singular value thresholding algorithm \cite{cai14} proposed to solve the problem RNuAMRM tends to lead to biased estimation by shrinking all the singular values toward to zero simultaneously, and sometimes results in over-penalization as the $\l_{1}$-norm in compressed sensing \cite{daubechies23}. Recently, some empirical evidence \cite{cui23}, has shown that, the non-convex algorithm, namely, iterative singular value thresholding (ISVT) algorithm, can really make a better recovery in some matrix rank minimization problems. However, the thresholding function for the ISVT algorithm is too complicated to computing, and converges slowly.

In this paper, inspired by the good performance of the $p$-thresholding operator \cite{voro24}, namely, "generalized thresholding operator" in compressed sensing, a generalized
singular value thresholding (GSVT) operator is generated to solve the problem ARMP. This GSVT operator comes to the soft thresholding operator \cite{cai14} when $p=1$, and for any $p<1$,
it penalizes small coefficients over a wider range and applies less bias to the larger coefficients, much like the hard thresholding operator \cite{Jian25} but without discontinuities.
With the change of the parameter $p$, we could get some much better results, which is one of the advantages for our algorithm compared with some state-of-art methods.

The rest of this paper is organized as the following. Some preliminaries knowledge that are used in this paper are given in Section \ref{section2}. Inspired by the generalized thresholding operator,
a generalized singular value thresholding operator is generated in Section \ref{section3}. In Section \ref{section4}, an iterative generalized singular value thresholding algorithm is proposed to solve the problem AMRM. In Section \ref{section5}, we demonstrate some numerical experiments on some matrix completion problems. Some conclusion remarks are presented in Section \ref{section6}.

\section{Preliminaries} \label{section2}
In this section, we give some preliminary knowledge that are used in this paper.

\subsection{Some notions} \label{subsection2-1}
For any matrix $X\in \mathbb{R}^{m\times n}$, let $X=U\Sigma_{X}V^{\top}=U[\mathrm{Diag}(\sigma(X)),\mathbf{0}]V^{\top}$ be the singular value decomposition (SVD) of matrix $X$, where $U$ is an $m\times m$ unitary matrix, $V$ is an $n\times n$ unitary matrix, $\Sigma_{X}=[\mathrm{Diag}(\sigma(X)),\mathbf{0}]\in \mathbb{R}^{m\times n}$ and $\sigma(X): \sigma_{1}(X)\geq\sigma_{2}(X)\geq\cdots\geq\sigma_{m}(X)$ denotes the singular value vector of matrix $X$, which is arranged in descending order. The linear map $\mathcal{A}$ determined
by $d$ given matrices $A_{1}, A_{2}, \cdots,A_{d}\in \mathbb{R}^{m\times n}$ is $\mathcal{A}(X)=(\langle A_{1},X\rangle, \langle A_{2},X\rangle,\cdots, \langle A_{d},X\rangle)^{\top}\in \mathbb{R}^{d}$. Let $\mathcal{A}^{\ast}$ denotes the adjoint of linear map $\mathcal{A}$. Then for any $y\in \mathbb{R}^{d}$, we have $\mathcal{A}^{\ast}(y)=\sum_{i=1}^{d}y_{i}A_{i}$. The standard inner product of matrices $X\in \mathbb{R}^{m\times n}$ and $Y\in \mathbb{R}^{m\times n}$ is given by
$\langle X,Y\rangle$, and $\langle X,Y\rangle=\mathrm{Tr}(Y^{\top}X)$. Define $A=(vec(A_{1}), vec(A_{2}),\cdots, vec(A_{d}))^{\top}\in \mathbb{R}^{d\times mn}$
and $x=vec(X)\in \mathbb{R}^{mn}$, we have $\mathcal{A}(X)=Ax$ and $\|\mathcal{A}(X)\|_{2}\leq\|\mathcal{A}\|_{2}\|X\|_{F}$.

\subsection{The form of AMRM and some relaxation forms} \label{subsection2-1}

The form of the affine matrix rank minimization (AMRM) problem in mathematics is given, i.e.,

\begin{equation}\label{equ1}
(\mathrm{AMRM})\ \ \ \min_{X\in \mathbb{R}^{m\times n}} \mbox{rank}(X)\ \ s.t. \ \  \mathcal{A}(X)=b,
\end{equation}
where $X\in \mathbb{R}^{m\times n}$, $b\in \mathbb{R}^{d}$ and $\mathcal{A}: \mathbb{R}^{m\times n}\mapsto \mathbb{R}^{d}$ is a linear map. A special case of AMRM is the matrix
completion (MC) problem:
\begin{equation}\label{equ2}
(\mathrm{MC})\ \ \ \min_{X\in \mathbb{R}^{m\times n}} \mbox{rank}(X)\ \ s.t. \ \  X_{ij}=M_{ij}
\end{equation}
for all $(i,j)\in\Omega$, where the only information available about $M\in \mathbb{R}^{m\times n}$ is a sampled set of entries $M_{i,j}, (i,j)\in\Omega$, and $\Omega$ is a
subset of the complete set of entries $\{1,2,\cdots,m\}\times\{1,2,\cdots,n\}$.

As the most popular convex relaxation, the nuclear-norm affine matrix rank minimization (NuAMRM) problem is given, i.e.,
\begin{equation}\label{equ3}\emph{}
(\mathrm{NuAMRM})\ \ \  \min_{X\in \mathbb{R}^{m\times n}} \ \|X\|_{\ast}\ \ s.t. \ \  \mathcal{A}(X)=b,
\end{equation}
where $\|X\|_{\ast}=\sum_{i=1}^{m}\sigma_{i}(X)$ is nuclear-norm of matrix $X$, and $\sigma_{i}(X)$ presents the $i$-th largest singular value of matrix $X$ arranged in
descending order.

The regularization form for affine matrix rank minimization (RNuAMRM) problem is given by
\begin{equation}\label{equ4}
(\mathrm{RNuAMRM})\ \ \ \  \min_{X\in \mathbb{R}^{m\times n}} \Big\{\frac{1}{2}\|\mathcal{A}(X)-b\|_{2}^{2}+\lambda \|X\|_{\ast}\Big\}
\end{equation}
where $\lambda>0$ is the regularization parameter. In \cite{cai14}, a singular value thresholding operator, namely, soft thresholding operator
$$\mathcal{D}_{\lambda}(X)=UD_{\lambda}(\Sigma)V^{\top},\ \mathcal{D}_{\lambda}(\Sigma)=[\mathrm{Diag}((\sigma(X)-\lambda)_{+}),\mathbf{0}]$$
is introduced to solve the problem RNuAMRM, where $t_{+}$ is the positive part of $t$, and $t_{+}=\max\{0,t\}$.

Recently, Cui et al.\cite{cui23} substituted the rank function $\mbox{rank}(X)$ by a sum of the non-convex fraction functions
\begin{equation}\label{equ5}
P_{a}(X)=\sum_{i=1}^{m}\rho_{a}(\sigma_{i}(X))
\end{equation}
in terms of the singular values of matrix $X\in \mathbb{R}^{m\times n}$, and the non-convex function
$$\rho_{a}(t)=\frac{a|t|}{a|t|+1}, \ \ a>0$$
is the fraction function. Then, they translated the problem (AMRM) into a transformed AMRM (TrAMRM) which has the following form
\begin{equation}\label{equ6}
(\mathrm{TrAMRM})\ \ \  \min_{X\in \mathbb{R}^{m\times n}} \ P_{a}(X)\ \ s.t. \ \  \mathcal{A}(X)=b
\end{equation}
for the constrained problem and
\begin{equation}\label{equ7}
(\mathrm{RTrAMRM})\ \ \  \min_{X\in \mathbb{R}^{m\times n}} \Big\{\frac{1}{2}\|\mathcal{A}(X)-b\|_{2}^{2}+\lambda P_{a}(X)\Big\}
\end{equation}
for the regularization problem. Moreover, an iterative singular value thresholding (ISVT) algorithm is proposed to solve the problem (RTrAMRM).

\section{Generalized singular value thresholding operator} \label{section3}

Inspired by the good performances of the generalized thresholding operator \cite{voro24,char31} in compressed sensing and differ from the former
thresholding operators (see \cite{cai14,cui23}), in this section, a generalized singular value thresholding operator is generated to solve the problem AMRM.

\begin{definition}\label{def1}{\rm(see \cite{voro24,char31})}
For any $\lambda>0$, $p\leq1$ and $w_{i}\in \mathbb{R}$, the generalized thresholding operator is given by
\begin{equation}\label{equ8}
\mathcal{r}_{\lambda,p}(w_{i})=\mathrm{sign}(w_{i})\max\Big\{0, |w_{i}|-\lambda|w_{i}|^{p-1}\Big\}.
\end{equation}
\end{definition}

\begin{definition}\label{def2}
(Vector generalized thresholding operator) For any $\lambda>0$ and $w=(w_{1}, w_{2}, \cdots, w_{m})^{\top}\in \mathbb{R}^{m}$, the vector generalized thresholding operator
$\mathcal{R}_{\lambda,p}$ is defined as
\begin{equation}\label{equ9}
\mathcal{R}_{\lambda,p}(w)=(\mathcal{r}_{\lambda,p}(w_{1}), \mathcal{r}_{\lambda,p}(w_{2}), \cdots, \mathcal{r}_{\lambda,p}(w_{m}))^{\top},
\end{equation}
where $\mathcal{r}_{\lambda,p}$ is defined in Definition \ref{def1}.
\end{definition}

\begin{lemma}\label{lem1}{\rm(see \cite{char31})}
Suppose $\mathcal{r}_{\lambda,p}: [0,\infty)\rightarrow \mathbb{R}_{+}$ is continuous, satisfies $w_{i}\leq\xi\Rightarrow \mathcal{r}_{\lambda,p}(w_{i})=0$ for $\xi\geq0$, is
strictly increasing on $[\xi,\infty)$, and $\mathcal{r}_{\lambda,p}(w_{i})\leq w_{i}$. Then the threshold operator $\mathcal{R}_{\lambda,p}(w)$ is the proximal mapping of the
penalty function $F_{p}(x)=\sum_{i=1}^{N}f_{p}(x_{i})$
$$\mathcal{R}_{\lambda,p}(w)\triangleq\arg\min_{x\in \mathbb{R}^{n}}\Big\{\frac{1}{2}\|x-w\|_{2}^{2}+\lambda F_{p}(x)\Big\},$$
where $f_{p}$ is even, strictly increasing and continuous on $[0,\infty)$, differentiable on $(0,\infty)$,  and non-differentiable at 0 if and only if $\xi>0$ (in which
case $\partial f_{p}(0)=[-1,1]$). If $w_{i}-\mathcal{r}_{\lambda,p}(w_{i})$ is non-increasing on $[\xi, \infty)$, then $f_{p}$ is concave on $[0,\infty)$ and $f_{p}$ satisfies
the triangle inequality.
\end{lemma}

\begin{figure}[h!]
 \centering
 \includegraphics[width=0.35\textwidth]{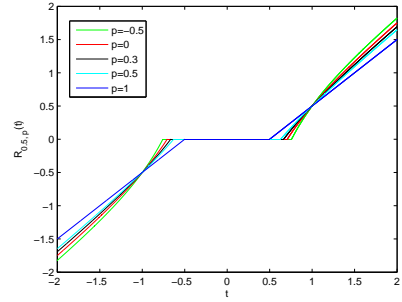}
\caption{Plot of generalized thresholding operator for some $p$ with $\lambda=0.5$.}
\label{figure1}
\end{figure}

To see clear that the generalized thresholding operator (\ref{equ8}) is equivalent to the classical soft thresholding operation \cite{dau28} in compressed sensing
when $p=1$, and it penalizes small coefficients over a wider range and applies less bias to the larger coefficients, much like the hard thresholding function \cite{blu29} but
without discontinuities for any $p<1$. The behavior of the generalized thresholding operator for some $p$ with $\lambda=0.5$ are painted in Figure \ref{figure1}.

\begin{remark}\label{remark1}
It is necessary to emphasize that the generalized thresholding operation defined in Definition \ref{def1} is an artificial operator and except $p=1$, it is
not the proximal mapping of the classical $\ell_{p}(0\leq p\leq1)$-norm minimization problem
\begin{equation}\label{equ10}
\min_{x\in \mathbb{R}^{n}}\Big\{\frac{1}{2}\|x-w\|_{2}^{2}+\lambda \|x\|_{p}^{p}\Big\}.
\end{equation}
The reasons can be detailed seen in \cite{xing30}.
\end{remark}

\begin{definition}\label{def3}
(Generalized singular value thresholding operator) For any $\lambda>0$, $p\leq1$ and $X=U[\mathrm{Diag}(\sigma(X)),\mathbf{0}]V^{\top}$ be the SVD of
matrix $X\in \mathbb{R}^{m\times n}$, the generalized singular value thresholding (GSVT) operator for matrix $X$ is defined by
\begin{equation}\label{equ11}
\mathcal{R}_{\lambda,p}(X)=U[\mathrm{Diag}(\mathcal{R}_{\lambda,p}(\sigma(X))),\mathbf{0}]V^{\top},
\end{equation}
where $\mathcal{R}_{\lambda,p}$ is defined in Definition \ref{def2}.
\end{definition}

The GSVT operator $\mathcal{R}_{\lambda, p}$ defined in Definition \ref{def3} is simply apply the vector generalized thresholding operator to the
singular values of a matrix, and effectively shrinks them towards zero. It is to see clear that the rank of the output matrix $\mathcal{R}_{\lambda,p}(X)$
is lower than the rank of the input matrix $X$.

Next, we will conclude the most important conclusion in this paper which underlies the algorithm to be proposed.

\begin{theorem}\label{the1}
For any $\lambda>0$, $p\leq1$, $X\in \mathbb{R}^{m\times n}$ and the penalty function $F_{p}(X)$ is in terms of the singular values of matrix $X$ and $F_{p}(X)=\sum_{i=1}^{m}f_{p}(\sigma_{i}(X))$,
then the GSVT operator $\mathcal{R}_{\lambda,p}(X)$ defined in Definition \ref{def2} is the proximal mapping of the penalty function $F_{p}(X)$:
\begin{equation}\label{equ12}
\mathcal{R}_{\lambda,p}(Y)\triangleq\arg\min_{X\in \mathbb{R}^{m\times n}}\Big\{\frac{1}{2}\|X-Y\|_{2}^{2}+\lambda F_{p}(X)\Big\},
\end{equation}
where $f_{p}$ is even,strictly increasing and continuous on $[0,\infty)$, differentiable on $(0,\infty)$, and non-differentiable at 0 if and only if $\xi>0$ (in which case
$\partial f_{p}(0)=[-1,1]$). If $\sigma_{i}(X)-\mathcal{r}_{\lambda,p}(\sigma_{i}(X))$ is non-increasing on $[\xi, \infty)$, then $f_{p}$ is concave on $[0,\infty)$ and $F_{p}$
satisfies the triangle inequality.
\end{theorem}

We will need the following technical lemma which is the key for proving Theorem \ref{the1}.

\begin{lemma}\label{lem2}
(von Neumann's trace inequality) For any matrices $X, Y\in \mathbb{R}^{m\times n}$ $(m\leq n)$, $\mathrm{Tr}(X^{\top}Y)\leq\sum_{i=1}^{m}\sigma_{i}(X)\sigma_{i}(Y)$, where $\sigma(X)$
and $\sigma(Y)$ are the singular value vector of matrices $X$ and $Y$ respectively. The equality holds if and only if there exists unitary matrices $U$ and $V$ that such
$X=U[\mathrm{Diag}(\sigma(X)),\mathbf{0}]V^{\top}$ and $Y=U[\mathrm{Diag}(\sigma(Y)),\mathbf{0}]V^{\top}$ as the singular value decompositions of the matrices $X$ and $Y$ simultaneously.
\end{lemma}

We now proceed to a proof of Theorem \ref{the1}.\\
\textbf{proof} {\rm(of Theorem \ref{the1})}
Since $\sigma(X): \sigma_{1}(X)\geq\sigma_{2}(X)\geq\cdots\geq\sigma_{m}(X)\geq0$ are the singular values of matrix $X$, the minimization problem
\begin{equation}\label{equ13}
\min_{X\in \mathbb{R}^{m\times n}}\Big\{\frac{1}{2}\|X-Y\|_{F}^{2}+\lambda F_{p}(X)\Big\}
\end{equation}
can be rewritten as
$$\min_{\sigma(X)}\Big\{\frac{1}{2}\|X-Y\|_{F}^{2}+\lambda \sum_{i=1}^{m}f_{p}(\sigma_{i}(X))\Big\}.$$
By using the trace inequality in Lemma \ref{lem2}, we have
\begin{eqnarray*}
&&\|X-Y\|_{F}^{2}\\
&&=\mathrm{Tr}(X^{\top}X)-2\mathrm{Tr}(X^{T}Y)+\mathrm{Tr}(Y^{\top}Y)\\
&&=\sum_{i=1}^{m}\sigma_{i}^{2}(X)-2\mathrm{Tr}(X^{\top}Y)+\sum_{i=1}^{m}\sigma_{i}^{2}(Y)\\
&&\geq\sum_{i=1}^{m}\sigma_{i}^{2}(X)-2\sum_{i=1}^{m}\sigma_{i}(X)\sigma_{i}(Y)+\sum_{i=1}^{m}\sigma_{i}^{2}(Y)\\
&&=\sum_{i=1}^{m}(\sigma_{i}(X)-\sigma_{i}(Y))^{2}.
\end{eqnarray*}
Noting that above equality holds if and only if the matrix $Y$ admits the singular value decomposition
$$Y=U[\mathrm{Diag}(\sigma(Y)),\mathbf{0}]V^{\top},$$
where $U$ and $V$ are the left and right orthonormal matrices in the SVD of matrix $X$ (see the last part of Lemma \ref{lem2}). So, the problem (\ref{equ13}) reduces to
\begin{equation}\label{equ14}
\min_{\sigma_{i}(X)}\Big\{\frac{1}{2}(\sigma_{i}(X)-\sigma_{i}(Y))^{2}+\lambda f_{p}(\sigma_{i}(X))\Big\},
\end{equation}
and following directly from the Lemma \ref{lem1} (Theorem 1 in \cite{char31}), we finish the proof.$\hfill{} \Box$

Although we do not know the exact expression of penalty function $F_{p}(X)$, the GSVT operator $\mathcal{R}_{\lambda,p}$ can really recover a low-rank matrix. In fact, the
GSVT operator equivalents to the soft thresholding operator \cite{cai14} when $p=1$. For $p<1$, it penalizes small coefficients over a wider range and applies less
bias to the larger coefficients, much like the hard thresholding operator \cite{Jian25} but without discontinuities. With the change of the parameter $p$, we may get
some much better results, which is one of the advantages for the GSVT operator compared with some other thresholding operator. This is the reason why we refer to
this transformation as the generalized thresholding operator for the GSVT operator.

\section{Iterative generalized singular value thresholding algorithm} \label{section4}

In this section, an iterative generalized singular value thresholding (IGSVT) algorithm is proposed to solve the problem AMRM. Moreover, the cross-validation
method \cite{xu33} is applied to adjust the regularization parameter $\lambda$ in each iteration.

\subsection{Fixed point inclusion for minimizer} \label{section4-1}

Now, we begin to consider our following regularization problem:
\begin{equation}\label{equ15}
\Psi_{1}(X)=\frac{1}{2}\|\mathcal{A}(X)-b\|_{2}^{2}+\lambda F_{p}(X).
\end{equation}
For any $\lambda>0$, $p<1$ and matrix $Z\in \mathbb{R}^{m\times n}$, let
\begin{equation}\label{equ16}
\begin{array}{llll}
\Psi_{2}(X,Z)&=&\mu\Big[\Psi_{1}(X)-\frac{1}{2}\|\mathcal{A}(X)-\mathcal{A}(Z)\|_{2}^{2}\Big]\\
&&+\frac{1}{2}\|X-Z\|_{F}^{2}.
\end{array}
\end{equation}
It is easy to verify that $\Psi_{2}(X, X)=\mu \Psi_{1}(X)$.

\begin{theorem}\label{the2}
For any positive numbers $\lambda>0$, $\mu>0$ and matrix $Z\in \mathbb{R}^{m\times n}$, if matrix $\tilde{X}\in \mathbb{R}^{m\times n}$ is the optimal solution of
the problem $\displaystyle\min_{X\in \mathbb{R}^{m\times n}}\Psi_{2}(X,Z)$, then
$$\tilde{X}= \mathcal{R}_{\lambda\mu, p}(B_{\mu}(Z))=\tilde{U}[\mathrm{Diag}(\mathcal{R}_{\lambda\mu,p}(\sigma(B_{\mu}(Z)))_{i}),\mathbf{0}]\tilde{V}^{\top},$$
where $B_{\mu}(Z)=Z-\mu \mathcal{A}^{\ast}\mathcal{A}(Z)+\mu \mathcal{A}^{\ast}(b)$, $B_{\mu}(Z)=\tilde{U}[\mathrm{Diag}(\sigma(B_{\mu}(Z))),\mathbf{0}]\tilde{V}^{\top}$
is the singular value decomposition of matrix $B_{\mu}(Z)$, the matrices $\tilde{U}$ and $\tilde{V}$ are the corresponding left and right orthonormal matrices, and
$\mathcal{R}_{\lambda\mu,p}$ is obtained by replacing $\lambda$ with $\lambda\mu$ in $\mathcal{R}_{\lambda,p}$.
\end{theorem}
\textbf{proof} By definition, the function $\Psi_{2}(X,Z)$ can be rewritten as
\begin{eqnarray*}
&&\Psi_{2}(X,Z)\\
&&=\frac{1}{2}\|X-(Z-\mu \mathcal{A}^{\ast}\mathcal{A}(Z)+\mu \mathcal{A}^{\ast}(b))\|_{F}^{2}\\
&&+\lambda\mu F_{p}(X)+\frac{\mu}{2}\|b\|_{2}^{2}+\frac{1}{2}\|Z\|_{F}^{2}-\frac{\mu}{2}\|\mathcal{A}(Z)\|_{2}^{2}\\
&&-\frac{1}{2}\|Z-\mu \mathcal{A}^{\ast}\mathcal{A}(Z)+\mu \mathcal{A}^{\ast}(b)\|_{F}^{2}\\
&&=\frac{1}{2}\|X-B_{\mu}(Z)\|_{F}^{2}+\lambda\mu F_{p}(X)+\frac{\mu}{2}\|b\|_{2}^{2}+\|Z\|_{F}^{2}\\
&&-\frac{\mu}{2}\|\mathcal{A}(Z)\|_{2}^{2}-\frac{1}{2}\|B_{\mu}(Z)\|_{F}^{2},
\end{eqnarray*}
which means that minimizing the function $\Psi_{2}(X,Z)$ on $X$, for any $\lambda>0$, $\mu>0$ and matrix $Z\in \mathbb{R}^{m\times n}$, is equivalent to
$$\min_{X\in \mathbb{R}^{m\times n}}\Big\{\frac{1}{2}\|X-B_{\mu}(Z)\|_{F}^{2}+\lambda\mu F_{p}(X)\Big\}.$$
By Theorem \ref{the1}, it is easy to verify that $\tilde{X}\in \mathbb{R}^{m\times n}$ is the optimal solution of the problem $\displaystyle\min_{X\in \mathbb{R}^{m\times n}}\Psi_{2}(X,Z)$
if and only if, for any $i$, $\sigma_{i}(\tilde{X})$ solves the problem
$$\min_{\sigma_{i}(X)}\Big\{\frac{1}{2}(\sigma_{i}(X)-\sigma_{i}(B_{\mu}(Z)))^{2}+\lambda f_{p}(\sigma_{i}(X))\Big\}.$$
Combing with Lemma \ref{lem1}, we finish this proof. $\hfill{} \Box$

Furthermore, if we take the parameter $\mu>0$ properly, we have
\begin{theorem}\label{the3}
For any positive numbers $\lambda>0$ and $0<\mu<\frac{1}{\|\mathcal{A}\|_{2}^{2}}$. If $X^{\ast}$ is the optimal solution of $\displaystyle\min_{X\in \mathbb{R}^{m\times n}}\Psi_{1}(X)$,
it can be expressed as
\begin{equation}\label{equ17}
\begin{array}{llll}
X^{\ast}&=&\mathcal{R}_{\lambda\mu, p}(B_{\mu}(X^{\ast}))\\
&=&U^{\ast}[\mathrm{Diag}(\mathcal{R}_{\lambda\mu,p}(\sigma(B_{\mu}(X^{\ast})))),\mathbf{0}](V^{\ast})^{\top},
\end{array}
\end{equation}
where $B_{\mu}(X^{\ast})=U^{\ast}[\mathrm{Diag}(\sigma(B_{\mu}(X^{\ast}))),\mathbf{0}](V^{\ast})^{\top}$ is the singular value decomposition of matrix $B_{\mu}(X^{\ast})$, the matrices $U^{\ast}$
and $V^{\ast}$ are the corresponding left and right orthonormal matrices, and $\mathcal{R}_{\lambda\mu,p}$ is obtained by replacing $\lambda$ with $\lambda\mu$ in $\mathcal{R}_{\lambda,p}$.
\end{theorem}
\textbf{proof} By definition of $\Psi_{2}(X, Y)$, we have
\begin{eqnarray*}
&&\Psi_{2}(X,X^{\ast})\\
&&=\mu\Big[\Psi_{1}(X)-\frac{1}{2}\|\mathcal{A}(X)-\mathcal{A}(X^{\ast})\|_{2}^{2}\Big]+\frac{1}{2}\|X-X^{\ast}\|_{F}^{2}\\
&&=\mu\Big[\frac{1}{2}\|\mathcal{A}(X)-b\|_{2}^{2}+\lambda F_{p}(X)\Big]+\frac{1}{2}\|X-X^{\ast}\|_{F}^{2}\\
&&-\frac{\mu}{2}\|\mathcal{A}(X)-\mathcal{A}(X^{\ast})\|_{2}^{2}\\
&&\geq\mu\Big[\frac{1}{2}\|\mathcal{A}(X)-b\|_{2}^{2}+\lambda F_{p}(X)\Big]\\
&&=\mu \Psi_{1}(X)\\
&&\geq\mu \Psi_{1}(X^{\ast})\\
&&=\Psi_{2}(X^{\ast},X^{\ast}),
\end{eqnarray*}
where the first inequality holds by the fact that
$$\|\mathcal{A}(X)-\mathcal{A}(X^{\ast})\|_{2}^{2}\leq\|\mathcal{A}\|_{2}^{2}\cdot\|X-X^{\ast}\|_{F}^{2}.$$
Combined with Theorem \ref{the2} and Theorem \ref{the1}, we can immediately finish this proof. $\hfill{} \Box$

Theorem \ref{the3} show us that, for any $0<\mu<\frac{1}{\|\mathcal{A}\|_{2}^{2}}$, if $X^{*}$ is the optimal solution of $\displaystyle\min_{X\in \mathbb{R}^{m\times n}}\Psi_{1}(X)$,
it also solves the problem $\displaystyle\min_{X\in \mathbb{R}^{m\times n}}\Psi_{2}(X,Z)$ with $Z=X^{\ast}$.

With the fixed point inclusion (\ref{equ17}), the IGSVT algorithm for solving the problem $\displaystyle\min_{X\in \mathbb{R}^{m\times n}}\Psi_{1}(X)$ can be naturally given by
\begin{equation}\label{equ18}
X^{k+1}=\displaystyle \mathcal{R}_{\lambda\mu, p}(B_{\mu}(X^{k})),\ \ \ \ k=0,1,\cdots,
\end{equation}
where $B_{\mu}(X^{k})=X^{k}-\mu \mathcal{A}^{\ast}\mathcal{A}(X^{k})+\mu \mathcal{A}^{\ast}(b)$.

\begin{algorithm}[h!]
\caption{: IGSVT algorithm}
\label{alg:A}
\begin{algorithmic}
\STATE {\textbf{input}: $\mathcal{A}: \mathbb{R}^{m\times n}\mapsto \mathbb{R}^{d}$, $b\in \mathbb{R}^{d}$}
\STATE {\textbf{initialize}: Given $X^{0}\in \mathbb{R}^{m\times n}$, $\mu=\frac{1-\varepsilon}{\|\mathcal{A}\|_{2}^{2}}(0<\varepsilon<1)$, $\lambda_{0}>0$ and $p\leq1$;}
\STATE {\textbf{while} not converged \textbf{do}}
\STATE \ \ \ \ \ \ \ {$Z^{k}:=B_{\mu}(X^{k})=X^{k}-\mu \mathcal{A}^{\ast}\mathcal{A}(X^{k})+\mu \mathcal{A}^{\ast}(b)$;}
\STATE {Compute the SVD of $Z^{k}$ as:}
\STATE \ \ \ \ \ \ \ {$Z^{k}:=U^{k}[\mathrm{Diag}(\sigma_{i}(Z^{k})),\mathbf{0}](V^{k})^{\top}$;}
\STATE \ \ \ \ \ \ \ {$\lambda=\lambda_{0}$;}
\STATE \ \ \ \ \ \ \ \ \ \ \ \ {\textbf{for}\ $i=1:m$}
\STATE \ \ \ \ \ \ \ \ \ \ \ \ \ \ \ {$\sigma_{i}(Z^{k+1})=\mathcal{r}_{\lambda\mu,p}(\sigma_{i}(Z^{k}))$}
\STATE \ \ \ \ \ \ \ {\ \ \ \ \ \textbf{end}}
\STATE \ \ \ \ \ \ \ {$X^{k+1}:=\mathcal{R}_{\lambda_{k}\mu, p}(Z^{k})=U^{k}[\mathrm{Diag}(\sigma(Z^{k+1})),\mathbf{0}](V^{k})^{\top}$}
\STATE \ \ \ \ \ \ \ {$k\rightarrow k+1$}
\STATE{\textbf{end while}}
\STATE{\textbf{return}: $X^{k+1}$}
\end{algorithmic}
\end{algorithm}

The following theorem establishes the convergence of IGSVT algorithm. Its proof follows from the specific condition that the step size $\mu$ satisfying $0<\mu<\frac{1}{\|\mathcal{A}\|_{2}^{2}}$
and a similar argument as used in the proof of [20, Theorem 3]

\begin{theorem}\label{the4}
Suppose the step size $\mu$ satisfying $0<\mu<\frac{1}{\|\mathcal{A}\|_{2}^{2}}$. Let the sequence $\{X^{k}\}$ be generated by IGSVT algorithm. There hold:
\begin{description}
  \item[$\mathrm{1)}$] The sequence $\{\Psi_{1}(X^{k})\}$ is decreasing.
  \item[$\mathrm{2)}$] $\{X^{k}\}$ is asymptotically regular, i.e., $\lim_{k\rightarrow\infty}\|X^{k+1}-X^{k}\|_{F}^{2}=0$.
  \item[$\mathrm{3)}$] Any accumulation point of $\{X^{k}\}$ is a stationary point.
\end{description}
\end{theorem}

\subsection{Adjusting values for the regularization parameter} \label{section4-2}

One problem needs to be addressed is that the IGSVT algorithm seriously depends on the setting of the regularization parameter $\lambda>0$. In this paper, the cross-validation
method \cite{xu33} is applied to adjust the regularization parameter $\lambda$ in each iteration. To make it clear, we suppose the matrix $X^{\ast}$ of rank $r$ is the
optimal solution to the problem $\displaystyle\min_{X\in \mathbb{R}^{m\times n}}\Psi_{1}(X)$. Define the singular
values of matrix $B_{\mu}(X^{\ast})$ as
$$\sigma_{1}(B_{\mu}(X^{\ast}))\geq\sigma_{2}(B_{\mu}(X^{\ast}))\geq\cdots\geq\sigma_{m}(B_{\mu}(X^{\ast})).$$
By equation (\ref{equ8}), we have
$$\sigma_{i}(B_{\mu}(X^{\ast}))>(\lambda\mu)^{\frac{1}{2-p}}\Leftrightarrow i\in\{1,2,\cdots,r\},$$
$$\sigma_{i}(B_{\mu}(X^{\ast}))\leq (\lambda\mu)^{\frac{1}{2-p}}\Leftrightarrow i\in\{r+1,r+2,\cdots,m\},$$
which implies
\begin{equation}\label{equ19}
\frac{(\sigma_{r+1}(B_{\mu}(X^{\ast})))^{2-p}}{\mu}\leq\lambda<\frac{(\sigma_{r}(B_{\mu}(X^{\ast})))^{2-p}}{\mu}.
\end{equation}
In practice, we approximate $\sigma_{i}((B_{\mu}(X^{\ast})))$ by $\sigma_{i}((B_{\mu}(X^{k})))$ in (\ref{equ19}),
and a choice of $\lambda$ is
\begin{equation}\label{equ20}
\lambda\in\bigg[\frac{(\sigma_{r+1}(B_{\mu}(X^{k})))^{2-p}}{\mu}, \frac{(\sigma_{r}(B_{\mu}(X^{k})))^{2-p}}{\mu}\bigg).
\end{equation}
Especially, we set
\begin{equation}\label{equ21}
\lambda=\lambda_{k}=\frac{(\sigma_{r+1}(B_{\mu}(X^{k})))^{2-p}}{\mu}
\end{equation}
in each iteration. That is, (\ref{equ21}) can be used to adjust the value of the regularization parameter $\lambda$ during iteration.

\section{Numerical experiments} \label{section5}
In the section, we first carry out a series of simulations to demonstrate the performances of the IGSVT algorithm on random low-rank matrix completion problems, and then
compared them with some other methods (singular value thresholding (SVT) algorithm \cite{cai14} and iterative singular value thresholding (ISVT) algorithm \cite{cui23}) on image
inpainting problems.

Two quantities are defined to quantify the difficulty of the low rank matrix recovery problems: $\mathrm{SR}=s/mn$ denotes the sampling ration,
where $s$ is the cardinality of observation set $\Omega$ whose entries are sampled randomly; $\mathrm{FR}=s/r(m+n-r)$ is the freedom ration, which is the ratio between
the number of sampled entries and the 'true dimensionality' of a $m\times n$ matrix of rank $r$, and it is a good quantity as the information oversampling ratio. In fact,
if $\mathrm{FR}<1$, it is impossible to recover an original low-rank matrix because there are an infinite number of matrices of rank $r$ with the observed entries \cite{ma32}.
The stopping criterion is usually as following
$$\frac{\|X^{k}-X^{k-1}\|_{F}}{\|X^{k}\|_{F}}\leq \mathrm{Tol},$$
where $X^{k}$ and $X^{k-1}$ are numerical results from two continuous iterative steps and $\mathrm{Tol}$ is a given small number. We set $\mathrm{Tol}=10^{-7}$ in our experiments.
In addition, the accuracy of the generated solution $X^{\ast}$ of our algorithm is measured by the relative error ($\mathrm{RE}$), which is defined as
$$\mathrm{RE}=\frac{\|X^{\ast}-M\|_{F}}{\|M\|_{F}},$$
where $M\in \mathbb{R}^{m\times n}$ is the given low-rank matrix.

\subsection{Completion of random matrices}
For the sake of simplicity, we set $m=n$ and generate $n\times n$ matrices $M$ of rank $r$ as the matrix products of two low-rank matrices $M_{1}$ and $M_{2}$ where
$M_{1}\in \mathbb{R}^{n\times r}$, $M_{2}\in \mathbb{R}^{r\times n}$ are generated with independent identically distributed Gaussian entries and the matrix $M=M_{1}M_{2}$
has rank at most $r$. To determine the best choice of parameter $p$ , we test IGSVT algorithm on random matrix completion problems with some different $p\in\{-0.9,-0.7, -0.5,-0.3,-0.2,-0.1,0.1,0.2,0.3,0.5,0.7,0.9\}$.

\begin{table}[h!]\footnotesize
\centering
\setlength{\tabcolsep}{0.5mm}{
\begin{tabular}{|c||l|l|l|l|l|l|}\hline
Problem&\multicolumn{2}{c}{$p=-0.9$}&\multicolumn{2}{|c}{$p=-0.7$}&\multicolumn{2}{|c|}{$p=-0.5$}\\
\hline
($n$,\,$r$,\,FR)&RE&Time&RE&Time&RE&Time\\
\hline
$(100,\,12,\,1.7730)$&1.01e-05& 11.37& 1.20e-05& 3.99& 8.24e-06& 2.07\\
\hline
$(200,\,12,\,3.4364)$&1.91e-06& 2.74& 2.43e-06& 2.01& 3.53e-06& 2.16\\
\hline
$(300,\,12,\,5.1020)$&1.54e-06& 3.15& 1.49e-06& 3.02& 1.42e-06& 2.92\\
\hline
$(400,\,12,\,6.7682)$&1.21e-06& 4.74& 1.21e-06& 4.90& 1.16e-06& 4.60\\
\hline
$(500,\,12,\,8.4345)$&9.88e-07& 6.31& 1.05e-06& 6.15& 1.19e-06& 6.93\\
\hline
$(600,\,12,\,10.1010)$&7.91e-07& 8.89& 8.45e-07& 8.91& 8.37e-07& 9.01\\
\hline
$(700,\,12,\,11.7675)$&7.70e-07& 12.30& 8.58e-07& 13.03& 8.65e-07& 12.12\\
\hline
$(800,\,12,\,13.4341)$&7.59e-07& 16.29& 7.85e-07& 16.39& 6.72e-07& 15.61\\
\hline
$(900,\,12,\,15.1007)$&6.83e-07& 21.15& 6.81e-07& 20.50& 7.15e-07& 20.35\\
\hline
$(1000,\,12,\,16.7673)$&6.54e-07& 27.17& 8.72e-07& 28.49& 8.78e-07& 28.09\\
\hline
$(1100,\,12,\,18.4339)$&9.74e-07& 44.72& 8.92e-07& 36.89& 7.50e-07& 35.76\\
\hline
$(1200,\,12,\,20.1005)$&6.89e-07& 46.39& 6.56e-07& 45.63& 7.39e-07& 45.89\\
\hline
\end{tabular}}
\caption{\scriptsize Numerical results of IGSVT algorithm for matrix completion problems with different $n$, FR and $p$ but fixed rank $r$, SR=0.40.}\label{table1}
\end{table}

\begin{table}[h!]\footnotesize
\centering
\setlength{\tabcolsep}{0.5mm}{
\begin{tabular}{|c||l|l|l|l|l|l|}\hline
Problem&\multicolumn{2}{c}{$p=-0.3$}&\multicolumn{2}{|c}{$p=-0.2$}&\multicolumn{2}{|c|}{$p=-0.1$}\\
\hline
($n$,\,$r$,\,FR)&RE&Time&RE&Time&RE&Time\\
\hline
$(100,\,12,\,1.7730)$&9.46e-06& 1.67& 8.42e-06& 1.83& 1.07e-05& 1.44\\
\hline
$(200,\,12,\,3.4364)$&2.97e-06& 2.04& 2.06e-06& 1.75& 2.69e-06& 1.87\\
\hline
$(300,\,12,\,5.1020)$&1.65e-06& 2.95& 1.52e-06& 2.74& 1.46e-06& 2.60\\
\hline
$(400,\,12,\,6.7682)$&1.01e-06& 4.31& 1.02e-06& 4.31& 1.22e-06& 4.41 \\
\hline
$(500,\,12,\,8.4345)$&8.78e-07& 6.18& 1.00e-06& 6.20& 1.16e-06& 6.48\\
\hline
$(600,\,12,\,10.1010)$&9.99e-07& 9.30& 8.53e-07& 9.09& 1.17e-06& 10.30\\
\hline
$(700,\,12,\,11.7675)$&9.34e-07& 12.77& 8.06e-07& 11.88& 7.91e-07& 12.12\\
\hline
$(800,\,12,\,13.4341)$&7.73e-07& 16.27& 1.03e-06& 16.56& 1.01e-06& 18.58\\
\hline
$(900,\,12,\,15.1007)$&7.97e-07& 20.90& 7.20e-07& 20.94& 7.48e-07& 21.93\\
\hline
$(1000,\,12,\,16.7673)$&6.97e-07& 27.14& 8.04e-07& 27.94& 6.92e-07& 27.15\\
\hline
$(1100,\,12,\,18.4339)$&5.91e-07& 34.83& 7.08e-07& 35.60& 6.69e-07& 35.20\\
\hline
$(1200,\,12,\,20.1005)$&8.09e-07& 46.23& 9.26e-07& 49.00& 8.00e-07& 46.33\\
\hline
\end{tabular}}
\caption{\scriptsize Numerical results of IGSVT algorithm for matrix completion problems with different $n$, FR and $p$ but fixed rank $r$, SR=0.40.}\label{table2}
\end{table}

\begin{table}[h!]\footnotesize
\centering
\setlength{\tabcolsep}{0.5mm}{
\begin{tabular}{|c||l|l|l|l|l|l|}\hline
Problem&\multicolumn{2}{c}{$p=0.1$}&\multicolumn{2}{|c}{$p=0.2$}&\multicolumn{2}{|c|}{$p=0.3$}\\
\hline
($n$,\,$r$,\,FR)&RE&Time&RE&Time&RE&Time\\
\hline
$(100,\,12,\,1.7730)$&7.14e-06& 2.04& 1.11e-05& 1.70& 9.81e-06& 1.61\\
\hline
$(200,\,12,\,3.4364)$&2.07e-06& 1.68& 3.47e-06& 1.88& 2.27e-06& 1.56\\
\hline
$(300,\,12,\,5.1020)$&1.90e-06& 2.93& 1.72e-06& 2.64& 1.70e-06& 2.66\\
\hline
$(400,\,12,\,6.7682)$&1.00e-06& 3.96& 1.16e-06& 4.11& 1.07e-06& 4.00\\
\hline
$(500,\,12,\,8.4345)$&9.92e-07& 5.87& 1.06e-06& 5.83& 1.34e-06& 5.84\\
\hline
$(600,\,12,\,10.1010)$&8.29e-07& 8.29& 9.29e-07& 8.30& 1.00e-06& 8.43\\
\hline
$(700,\,12,\,11.7675)$&9.10e-07& 11.74& 7.34e-07& 11.52& 8.27e-07& 11.41\\
\hline
$(800,\,12,\,13.4341)$&7.64e-07& 15.15& 8.52e-07& 15.02& 8.07e-07& 15.18\\
\hline
$(900,\,12,\,15.1007)$&6.96e-07& 19.59& 7.03e-07& 19.82& 8.82e-07& 20.06\\
\hline
$(1000,\,12,\,16.7673)$&6.87e-07& 26.06& 7.42e-07& 25.96& 6.83e-07& 25.80\\
\hline
$(1100,\,12,\,18.4339)$&6.96e-07& 33.60& 6.56e-07& 33.92& 6.72e-07& 33.39\\
\hline
$(1200,\,12,\,20.1005)$&7.68e-07& 44.51& 8.10e-07& 46.97& 6.45e-07& 44.26\\
\hline
\end{tabular}}
\caption{\scriptsize Numerical results of IGSVT algorithm for matrix completion problems with different $n$, FR and $p$ but fixed rank $r$, SR=0.40.}\label{table3}
\end{table}

\begin{table}[h!]\footnotesize
\centering
\setlength{\tabcolsep}{0.5mm}{
\begin{tabular}{|c||l|l|l|l|l|l|}\hline
Problem&\multicolumn{2}{c}{$p=0.5$}&\multicolumn{2}{|c}{$p=0.7$}&\multicolumn{2}{|c|}{$p=0.9$}\\
\hline
($n$,\,$r$,\,FR)&RE&Time&RE&Time&RE&Time\\
\hline
$(100,\,12,\,1.7730)$&9.82e-06& 1.44& 7.45e-06& 1.13& 9.69e-06& 1.73\\
\hline
$(200,\,12,\,3.4364)$&2.95e-06& 1.58& 2.49e-06& 1.58& 3.37e-06& 2.28\\
\hline
$(300,\,12,\,5.1020)$&1.56e-06& 2.71& 1.41e-06& 2.53& 1.99e-06& 3.21\\
\hline
$(400,\,12,\,6.7682)$&1.41e-06& 3.98& 1.22e-06& 4.00& 1.15e-06& 4.48\\
\hline
$(500,\,12,\,8.4345)$&1.07e-06& 5.96& 1.28e-06& 6.42& 1.27e-06& 7.02\\
\hline
$(600,\,12,\,10.1010)$&8.83e-07& 8.27& 1.25e-06& 9.00& 1.25e-06& 9.81\\
\hline
$(700,\,12,\,11.7675)$&8.54e-07& 11.39& 8.28e-07& 11.67& 9.46e-07& 13.17\\
\hline
$(800,\,12,\,13.4341)$&8.05e-07& 14.86& 8.46e-07& 15.36& 9.11e-07& 16.82\\
\hline
$(900,\,12,\,15.1007)$&8.40e-07& 19.98& 8.34e-07& 20.10& 8.54e-07& 22.68\\
\hline
$(1000,\,12,\,16.7673)$&8.08e-07& 26.42& 6.93e-07& 26.27& 8.47e-07& 30.26\\
\hline
$(1100,\,12,\,18.4339)$&6.40e-07& 33.56& 5.84e-07& 33.92& 7.46e-07& 38.95\\
\hline
$(1200,\,12,\,20.1005)$& 5.99e-07& 44.34& 8.74e-07& 46.98& 6.65e-07& 48.92\\
\hline
\end{tabular}}
\caption{\scriptsize Numerical results of IGSVT algorithm for matrix completion problems with different $n$, FR and $p$ but fixed rank $r$, SR=0.40.}\label{table4}
\end{table}

\begin{table}[h!]\footnotesize
\centering
\setlength{\tabcolsep}{0.7mm}{
\begin{tabular}{|c||l|l|l|l|l|l|}\hline
Problem&\multicolumn{2}{c}{$p=-0.9$}&\multicolumn{2}{|c}{$p=-0.7$}&\multicolumn{2}{|c|}{$p=-0.5$}\\
\hline
($n$,\,$\tilde{r}$,\,FR)&RE&Time&RE&Time&RE&Time\\
\hline
$(100,\,11,\,1.9240)$&5.79e-06& 2.19& 5.90e-06& 2.37& 6.39e-06& 1.97\\
\hline
$(100,\,12,\,1.7730)$&7.76e-06& 3.82& 6.85e-06& 6.12& 1.21e-05& 2.30\\
\hline
$(100,\,13,\,1.6454)$&1.39e-05& 13.13& 1.25e-05& 5.64& 1.19e-05& 4.30\\
\hline
$(100,\,14,\,1.5361)$&1.99e-05& 13.24& 1.14e-05& 7.06& 1.75e-05& 11.46\\
\hline
$(100,\,15,\,1.4414)$&1.56e-05& 40.47& 1.80e-05& 13.08& 1.99e-05& 12.55\\
\hline
$(100,\,16,\,1.3587)$&---& ---& 2.61e-05& 13.17& 2.68e-05& 13.49\\
\hline
$(100,\,17,\,1.2858)$&---& ---& ---& ---& ---& ---\\
\hline
$(100,\,18,\,1.2210)$&---& ---& ---& ---& ---& ---\\
\hline
$(100,\,19,\,1.1631)$&---& ---& ---& ---& ---& ---\\
\hline
$(100,\,20,\,1.1111)$&---& ---& ---& ---& ---& ---\\
\hline
$(100,\,21,\,1.0641)$&---& ---& ---& ---& ---& ---\\
\hline
$(100,\,22,\,1.0215)$&---& ---& ---& ---& ---& ---\\
\hline
\end{tabular}}
\caption{\scriptsize Numerical results of IGSVT algorithm  for matrix completion problems with different rank $r$, FR and $p$ but fixed $n$, SR=0.40.}\label{table5}
\end{table}

\begin{table}[h!]\footnotesize
\centering
\setlength{\tabcolsep}{0.7mm}{
\begin{tabular}{|c||l|l|l|l|l|l|}\hline
Problem&\multicolumn{2}{c}{$p=-0.3$}&\multicolumn{2}{|c}{$p=-0.2$}&\multicolumn{2}{|c|}{$p=-0.1$}\\
\hline
($n$,\,$\tilde{r}$,\,FR)&RE&Time&RE&Time&RE&Time\\
\hline
$(100,\,11,\,1.9240)$&5.62e-06& 1.94& 6.19e-06& 2.36& 1.06e-05& 1.55\\
\hline
$(100,\,12,\,1.7730)$&9.18e-06& 1.78& 9.01e-06& 1.56& 8.88e-06& 1.58\\
\hline
$(100,\,13,\,1.6454)$&1.33e-05& 5.50& 3.23e-05& 4.03& 1.28e-05& 1.91\\
\hline
$(100,\,14,\,1.5361)$&1.72e-05& 4.01& 1.29e-05& 2.32& 1.32e-05& 1.97\\
\hline
$(100,\,15,\,1.4414)$&1.43e-05& 10.81& 1.95e-05& 3.38& 1.46e-05& 2.98\\
\hline
$(100,\,16,\,1.3587)$&4.68e-05& 7.81& 6.60e-05& 6.18& 1.90e-05& 3.87\\
\hline
$(100,\,17,\,1.2858)$&5.12e-05& 9.68& 3.48e-05& 7.29& 3.00e-05& 5.39\\
\hline
$(100,\,18,\,1.2210)$&---& ---& 5.29e-05& 23.60& 6.88e-05& 10.50\\
\hline
$(100,\,19,\,1.1631)$&---& ---& ---& ---& ---& ---\\
\hline
$(100,\,20,\,1.1111)$&---& ---& ---& ---& ---& ---\\
\hline
$(100,\,21,\,1.0641)$&---& ---& ---& ---& ---& ---\\
\hline
$(100,\,22,\,1.0215)$&---& ---& ---& ---& ---& ---\\
\hline
\end{tabular}}
\caption{\scriptsize Numerical results of IGSVT algorithm  for matrix completion problems with different rank $r$, FR and $p$ but fixed $n$, SR=0.40.}\label{table6}
\end{table}

\begin{table}[h!]\footnotesize
\centering
\setlength{\tabcolsep}{0.7mm}{
\begin{tabular}{|c||l|l|l|l|l|l|}\hline
Problem&\multicolumn{2}{c}{$p=0.1$}&\multicolumn{2}{|c}{$p=0.2$}&\multicolumn{2}{|c|}{$p=0.3$}\\
\hline
($n$,\,$\tilde{r}$,\,FR)&RE&Time&RE&Time&RE&Time\\
\hline
$(100,\,11,\,1.9240)$&8.17e-06& 1.12& 9.00e-06& 1.20& 1.25e-05& 1.47\\
\hline
$(100,\,12,\,1.7730)$&1.03e-05& 1.34& 1.52e-05& 1.60& 1.16e-05& 1.83\\
\hline
$(100,\,13,\,1.6454)$&1.44e-05& 1.61& 1.61e-05& 2.06& 1.14e-05& 1.69\\
\hline
$(100,\,14,\,1.5361)$&2.19e-05& 2.51& 1.33e-05& 1.85& 1.77e-05& 2.07\\
\hline
$(100,\,15,\,1.4414)$&1.49e-05& 2.30& 1.49e-05& 2.05& 2.05e-05& 2.21\\
\hline
$(100,\,16,\,1.3587)$&2.28e-05& 3.00& 3.01e-05& 2.84& 3.10e-05& 3.71\\
\hline
$(100,\,17,\,1.2858)$&2.58e-05& 3.83& 6.65e-05& 5.93& 9.41e-05& 3.57\\
\hline
$(100,\,18,\,1.2210)$&4.34e-05& 5.81& 5.30e-05& 5.53& 4.68e-05& 6.67\\
\hline
$(100,\,19,\,1.1631)$&7.61e-05& 11.10& 7.84e-05& 7.91& 9.27e-05& 12.73\\
\hline
$(100,\,20,\,1.1111)$&---& ---& 1.88e-04& 23.77& 1.38e-04& 11.75\\
\hline
$(100,\,21,\,1.0641)$&---& ---& ---& ---& ---& ---\\
\hline
$(100,\,22,\,1.0215)$&---& ---& ---& ---& ---& ---\\
\hline
\end{tabular}}
\caption{\scriptsize Numerical results of IGSVT algorithm  for matrix completion problems with different rank $r$, FR and $p$ but fixed $n$, SR=0.40.}\label{table7}
\end{table}

\begin{table}[h!]\footnotesize
\centering
\setlength{\tabcolsep}{0.7mm}{
\begin{tabular}{|c||l|l|l|l|l|l|}\hline
Problem&\multicolumn{2}{c}{$p=0.5$}&\multicolumn{2}{|c}{$p=0.7$}&\multicolumn{2}{|c|}{$p=0.9$}\\
\hline
($n$,\,$\tilde{r}$,\,FR)&RE&Time&RE&Time&RE&Time\\
\hline
$(100,\,11,\,1.9240)$&9.79e-06& 1.16& 9.67e-06& 1.22& 7.47e-06& 1.47\\
\hline
$(100,\,12,\,1.7730)$&9.95e-06& 1.53& 8.50e-06& 1.25& 1.14e-05& 2.17\\
\hline
$(100,\,13,\,1.6454)$&1.51e-05& 1.61& 1.85e-05& 171& 1.44e-05& 2.79\\
\hline
$(100,\,14,\,1.5361)$&1.45e-05& 1.65& 2.18e-05& 2.46& 1.78e-05& 4.84\\
\hline
$(100,\,15,\,1.4414)$&4.02e-05& 2.77& 1.84e-05& 2.36& 1.20e-02& 3.28\\
\hline
$(100,\,16,\,1.3587)$&2.94e-05& 2.67& 2.33e-05& 2.63& 2.06e-02& 1.59\\
\hline
$(100,\,17,\,1.2858)$&5.75e-05& 4.05& 3.10e-05& 3.37& 2.81e-02& 1.22\\
\hline
$(100,\,18,\,1.2210)$&4.50e-05& 4.10& 4.36e-05& 4.72& 2.78e-02& 1.19\\
\hline
$(100,\,19,\,1.1631)$&9.04e-05& 6.50& 9.62e-05& 7.25& 3.25e-02& 0.97\\
\hline
$(100,\,20,\,1.1111)$&1.25e-04& 9.38& 2.47e-04& 13.61& 3.32e-02& 1.06\\
\hline
$(100,\,21,\,1.0641)$&3.93e-04& 18.83& 4.08e-04& 34.17& 3.36e-02& 0.97\\
\hline
$(100,\,22,\,1.0215)$&2.10e-03& 38.02& 8.60e-03& 25.60& 3.54e-02& 1.09\\
\hline
\end{tabular}}
\caption{\scriptsize Numerical results of IGSVT algorithm  for matrix completion problems with different rank $r$, FR and $p$ but fixed $n$, SR=0.40.}\label{table8}
\end{table}

Tables \ref{table1}, \ref{table2}, \ref{table3}, \ref{table4} report the numerical results of IGSVT algorithm for the random low-rank matrix completion problems with $\mathrm{SR}=0.40$ when we fix rank $r=12$ and vary $n$ from $100$ to $1200$ with step size $100$. Tables \ref{table5}, \ref{table6}, \ref{table7}, \ref{table8} present numerical results of the IGSVT algorithms in the case where $n$ is fixed to 100 and $r$ is varied from 11 to 22 with step size 1. Comparing the performances of IGSVT algorithm for completion of random low rank matrices with different $p$ and $\mathrm{FR}$ we find that $p=0.5$ is the best strategy when $\mathrm{FR}$ is closed to one.

\subsection{Image inpainting}
\begin{figure}[h!]
  \centering
  \begin{minipage}[t]{0.45\linewidth}
  \centering
  \includegraphics[width=1\textwidth]{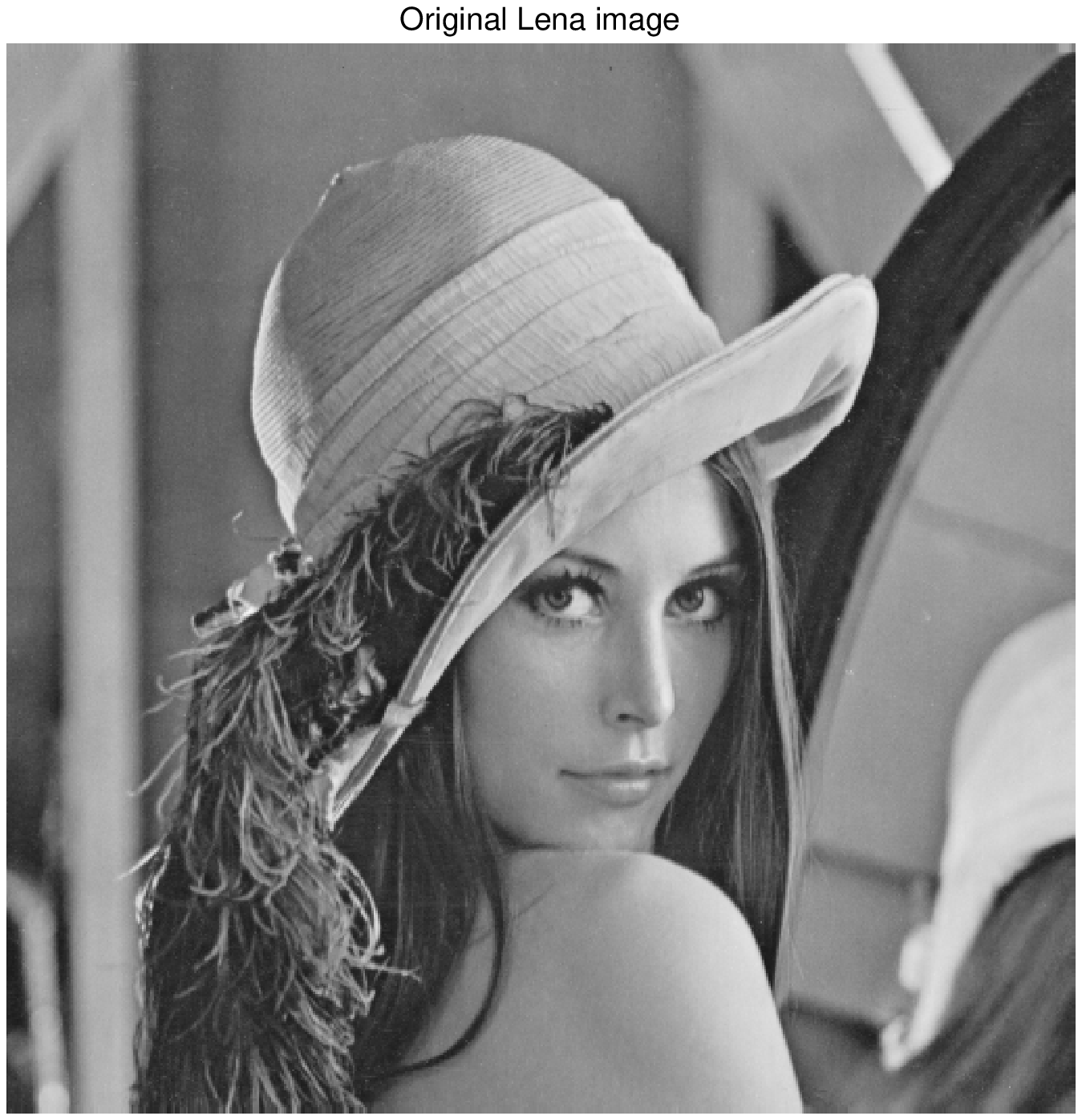}
  \end{minipage}
  \begin{minipage}[t]{0.45\linewidth}
  \centering
  \includegraphics[width=1\textwidth]{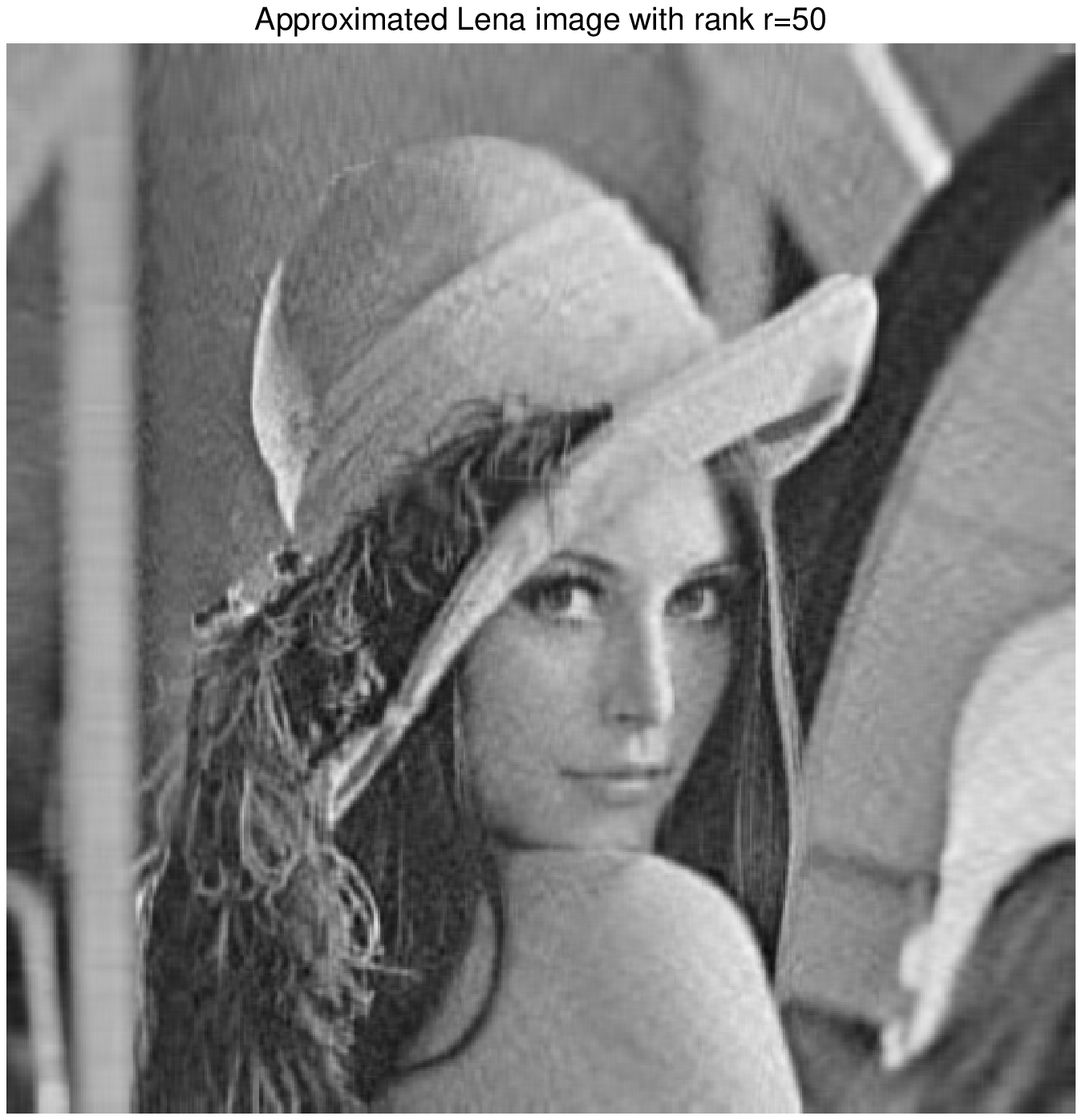}
  \end{minipage}
  \caption{Original $512\times512$ gray-scale Lena image and its approximation with rank $r=50$.}\label{figure2}
\end{figure}

\begin{figure}[h!]
  \centering
  \begin{minipage}[t]{0.45\linewidth}
  \centering
  \includegraphics[width=1\textwidth]{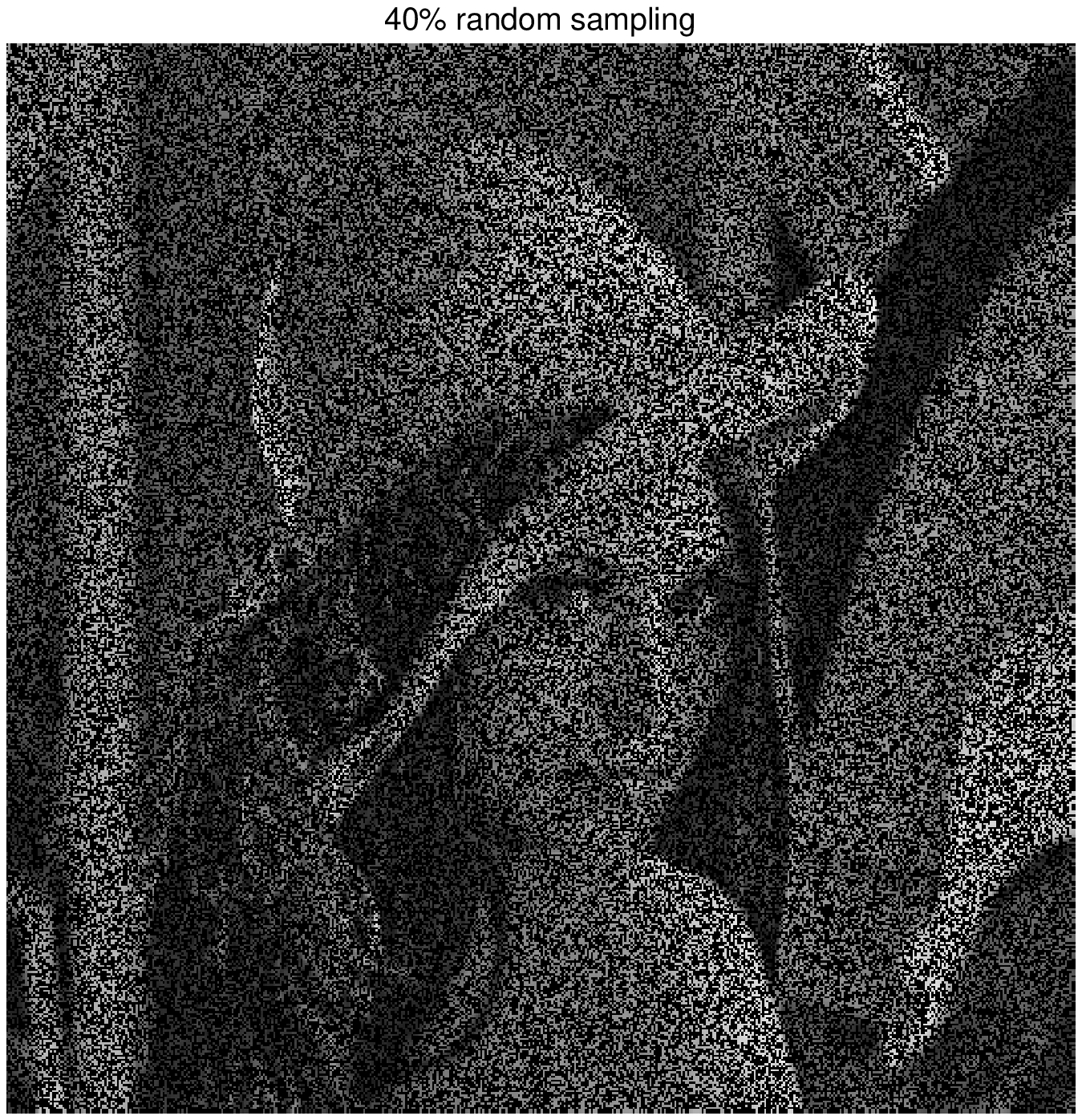}
  \end{minipage}
  \begin{minipage}[t]{0.45\linewidth}
  \centering
  \includegraphics[width=1\textwidth]{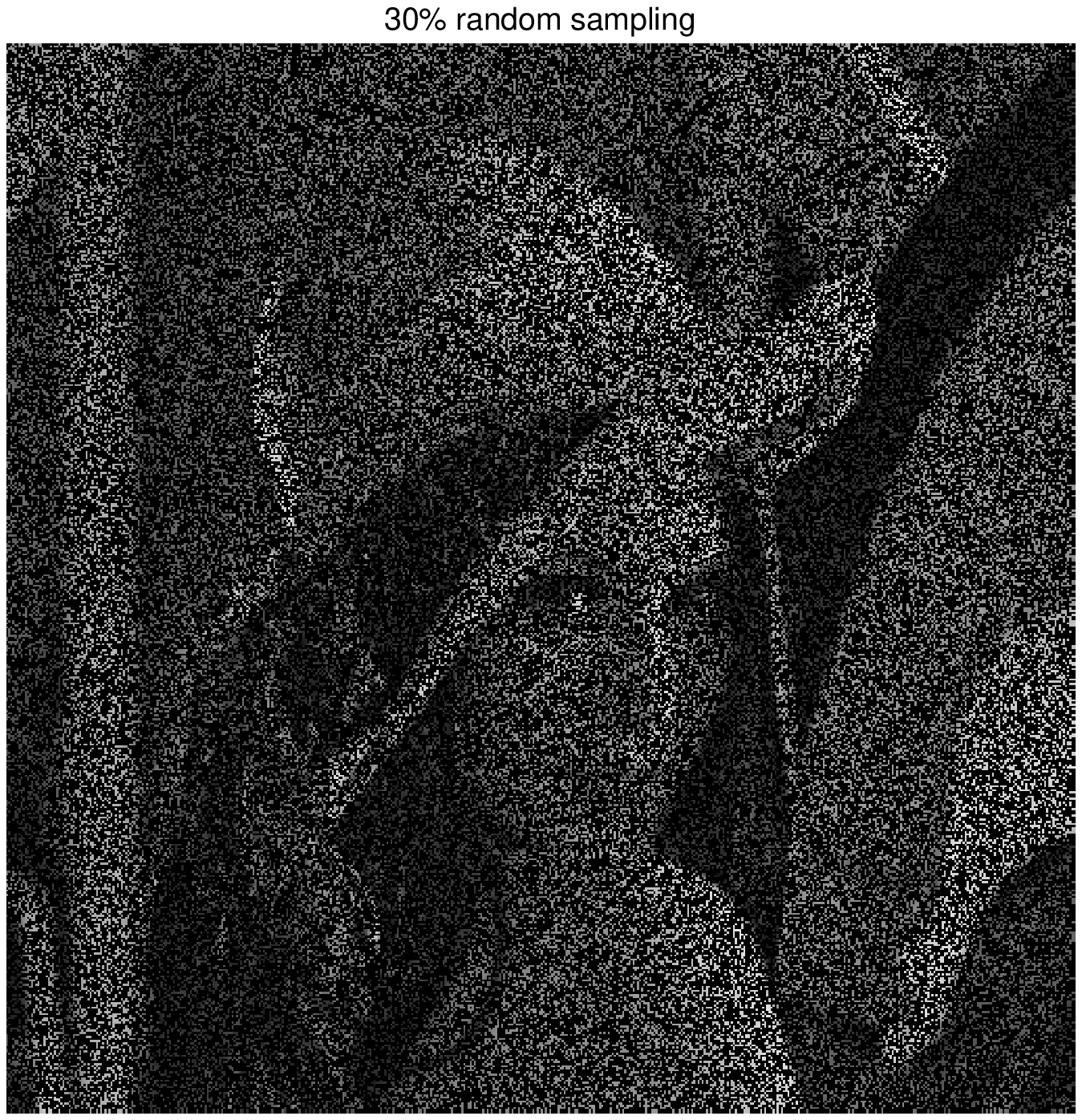}
  \end{minipage}
  \caption{Two sampled low-rank images (Left: $\mathrm{SR}=0.40$. Right: $\mathrm{SR}=0.30$).} \label{figure3}
\end{figure}

\begin{figure}[h!]
  \centering
  \begin{minipage}[t]{0.45\linewidth}
  \centering
  \includegraphics[width=1\textwidth]{lenasample04.eps}
  \end{minipage}
  \begin{minipage}[t]{0.45\linewidth}
  \centering
  \includegraphics[width=1\textwidth]{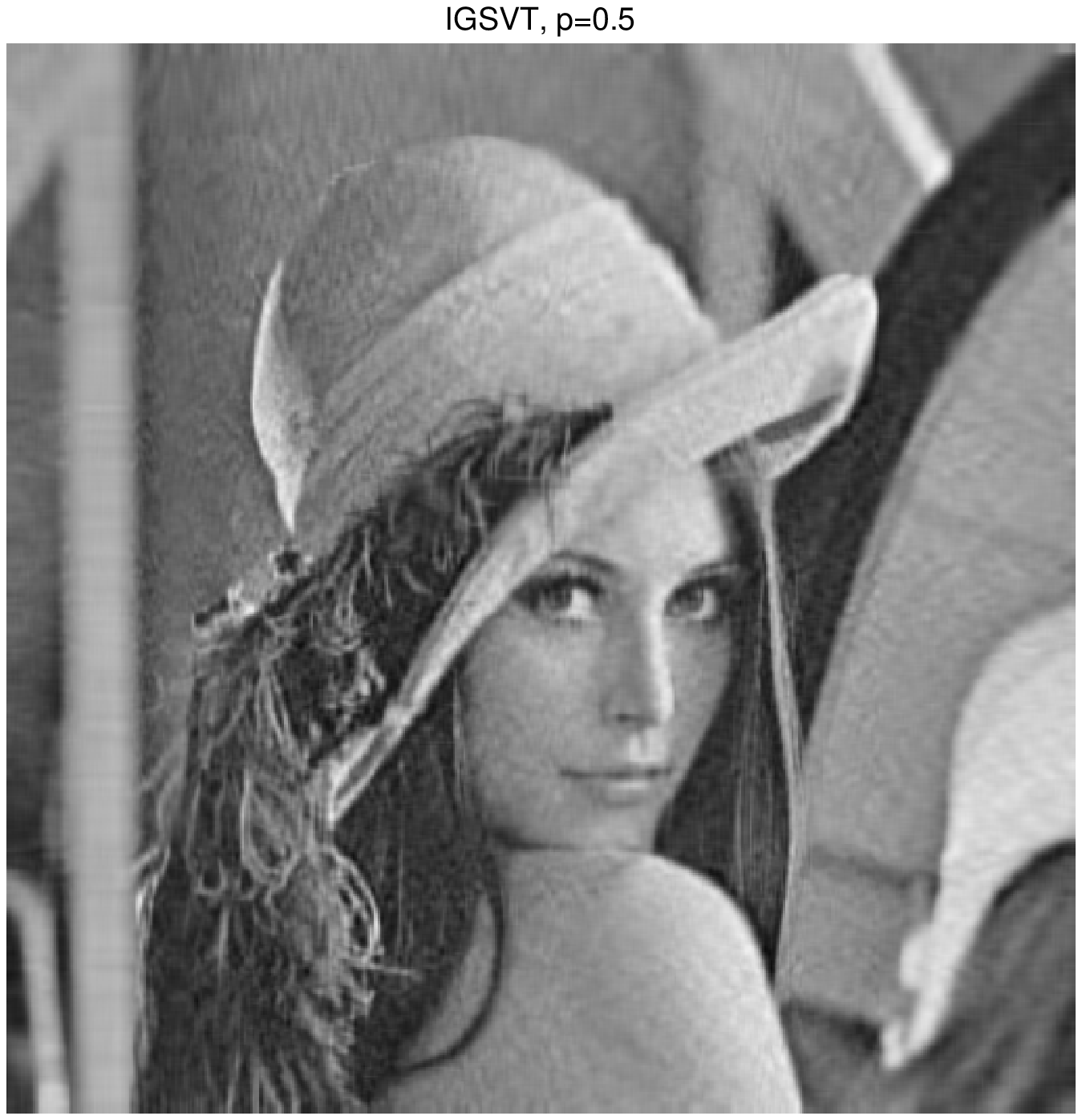}
  \end{minipage}
    \begin{minipage}[t]{0.45\linewidth}
  \centering
  \includegraphics[width=1\textwidth]{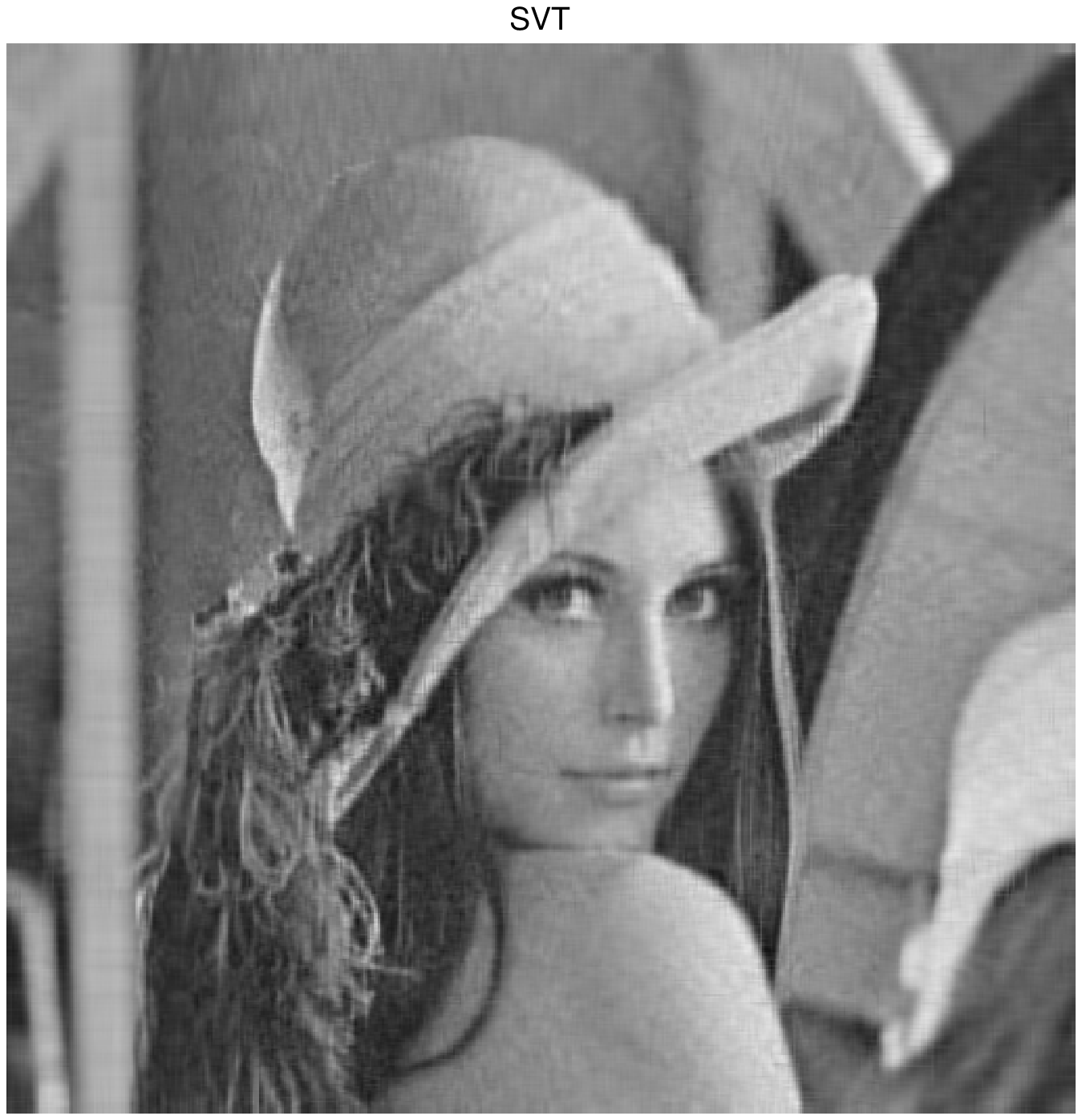}
  \end{minipage}
    \begin{minipage}[t]{0.45\linewidth}
  \centering
  \includegraphics[width=1\textwidth]{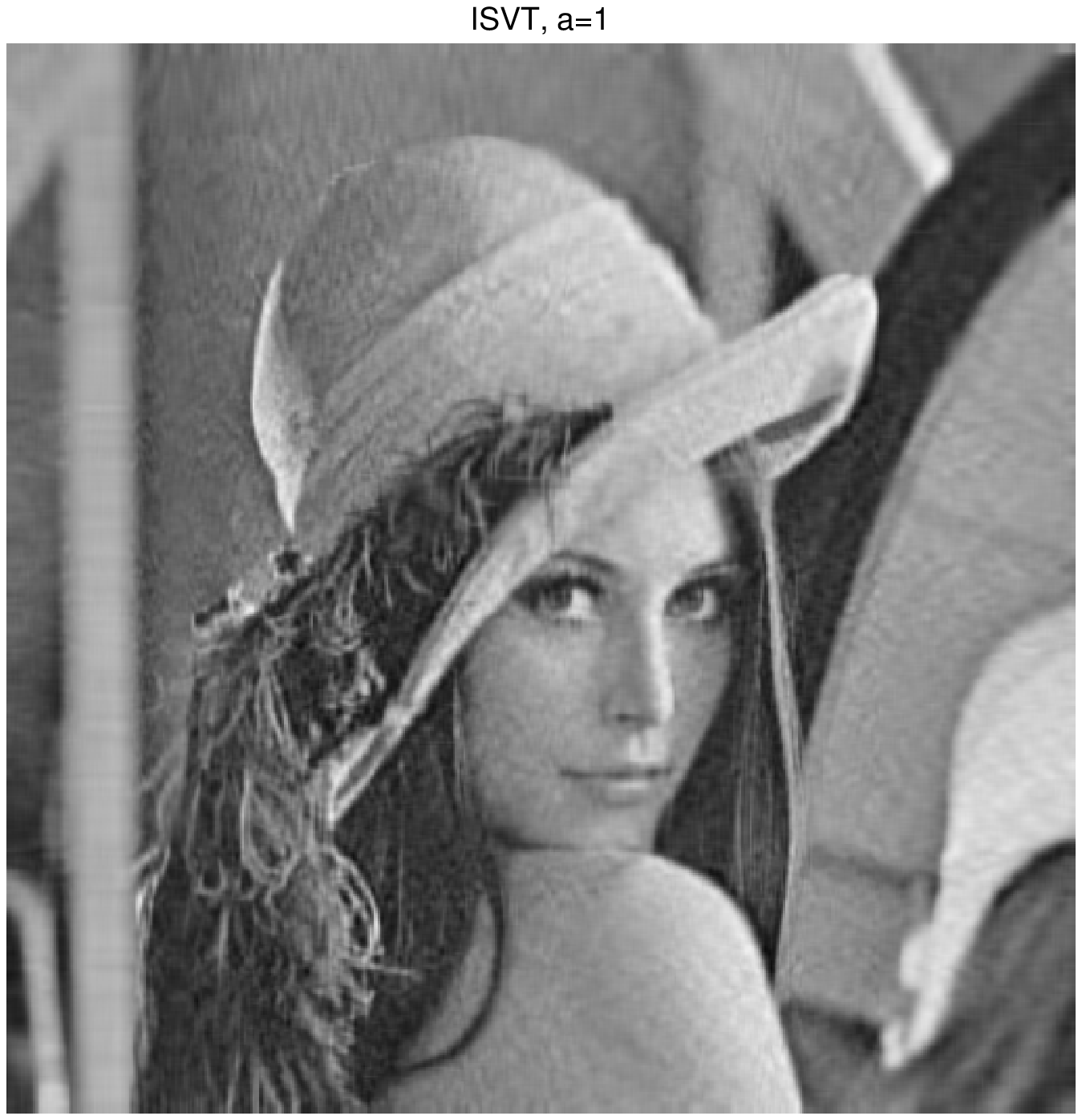}
  \end{minipage}
  \caption{Comparisons of IGSVT algorithm, SVT algorithm and ISVT algorithm for image inpainting with $\mathrm{SR}=0.40$.} \label{figure4}
\end{figure}

\begin{figure}[h!]
  \centering
  \begin{minipage}[t]{0.45\linewidth}
  \centering
  \includegraphics[width=1\textwidth]{lenasample03.eps}
  \end{minipage}
  \begin{minipage}[t]{0.45\linewidth}
  \centering
  \includegraphics[width=1\textwidth]{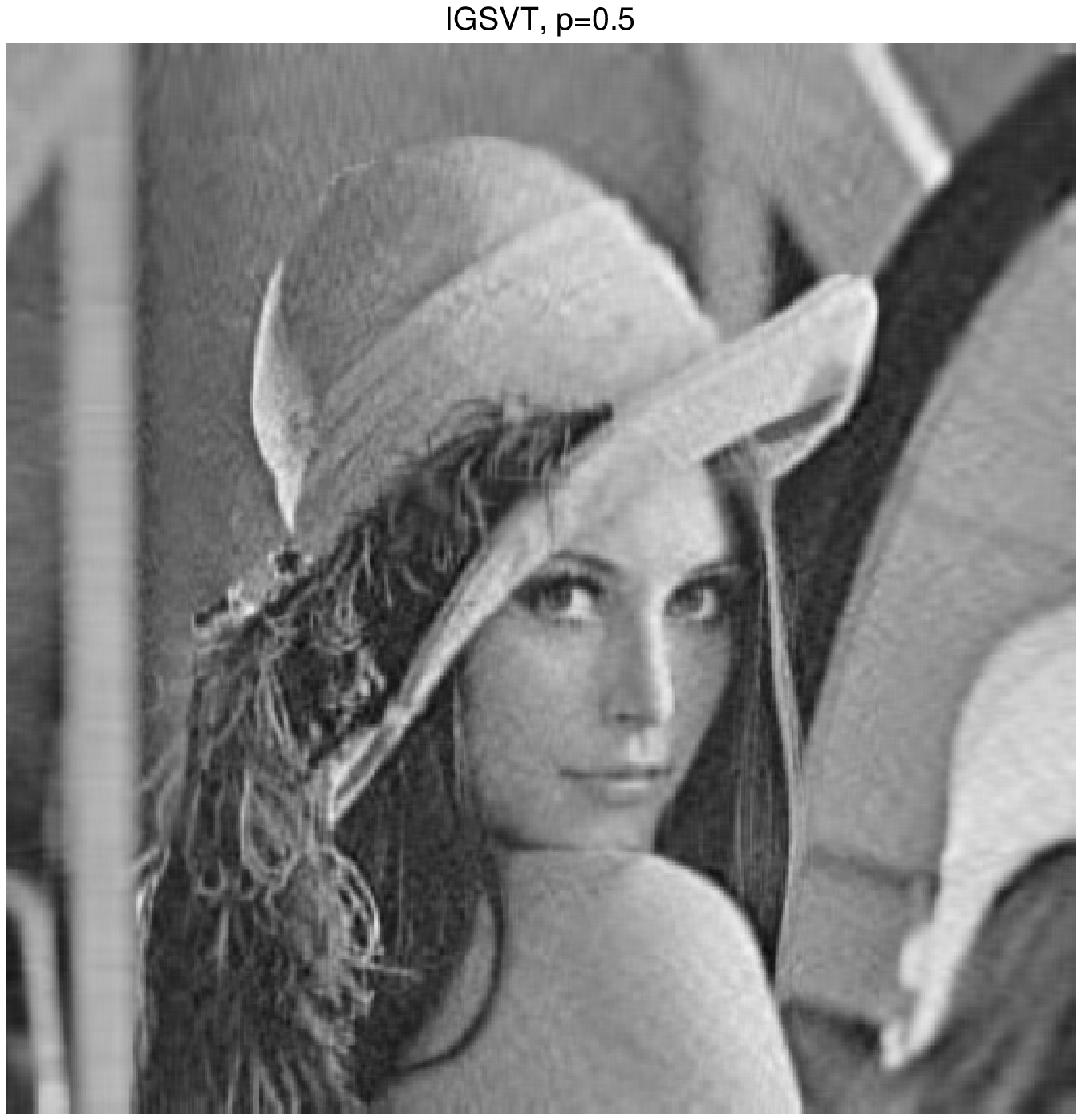}
  \end{minipage}
    \begin{minipage}[t]{0.45\linewidth}
  \centering
  \includegraphics[width=1\textwidth]{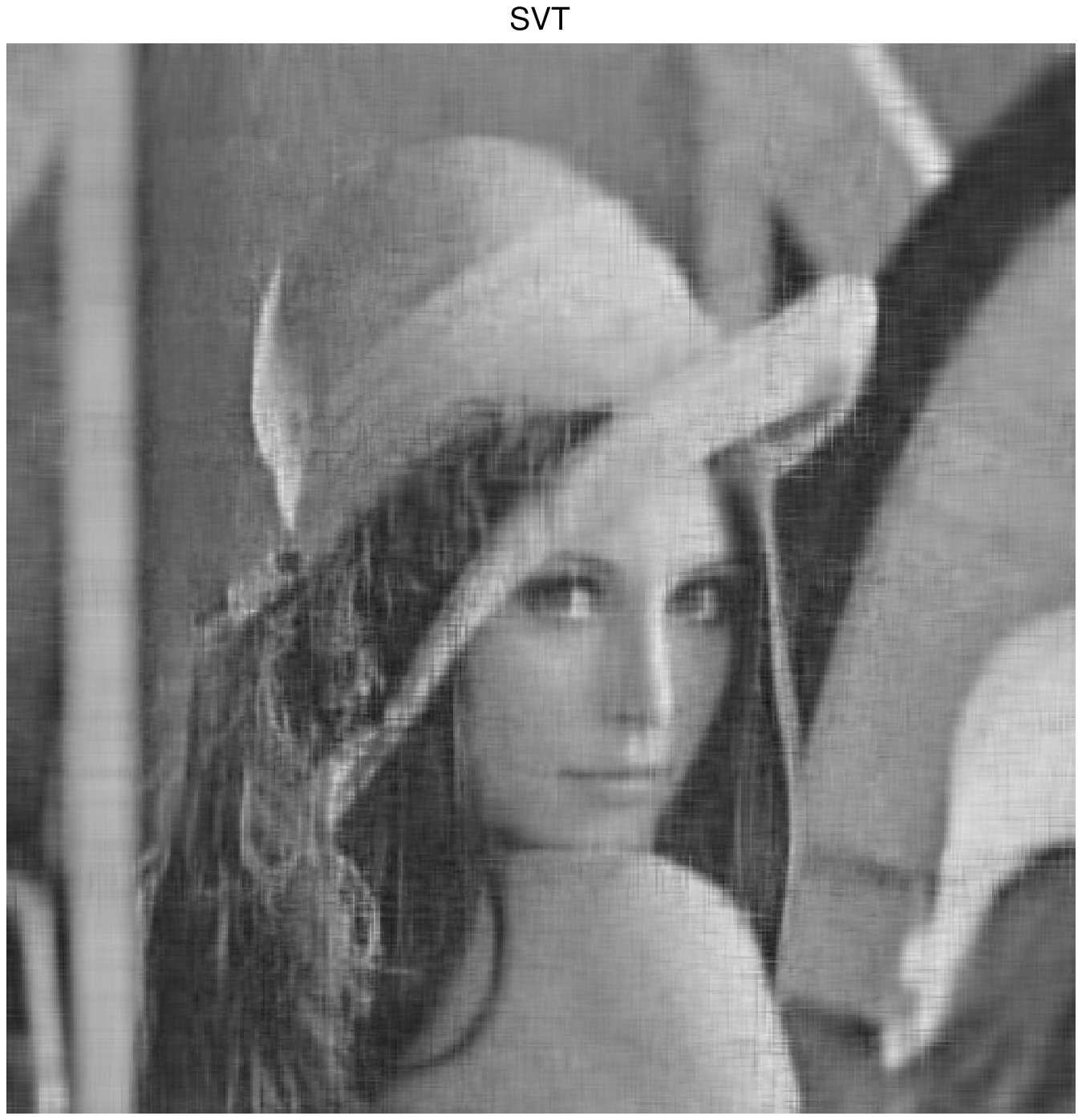}
  \end{minipage}
    \begin{minipage}[t]{0.45\linewidth}
  \centering
  \includegraphics[width=1\textwidth]{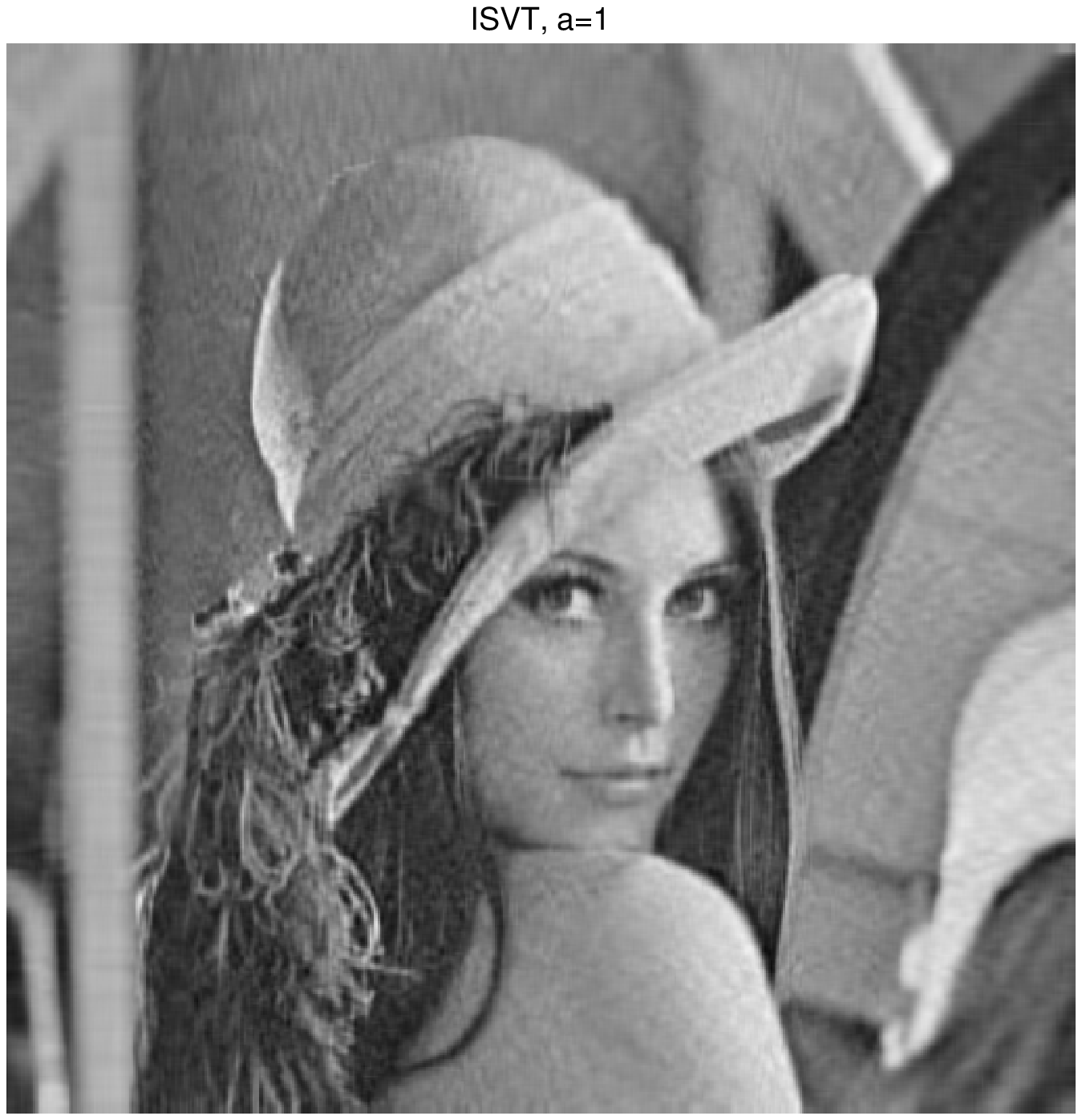}
  \end{minipage}
  \caption{Comparisons of IGSVT algorithm, SVT algorithm and ISVT algorithm for image inpainting with $\mathrm{SR}=0.30$.} \label{figure5}
\end{figure}

In the experiments, the IGSVT algorithm is tested on image inpainting problems and compared it with some state-of-art methods (singular value thresholding (SVT) algorithm \cite{cai14}
and iterative singular value thresholding (ISVT) algorithm \cite{cui23}). The three algorithms are tested on a standard $512\times512$ gray-scale image (Lena). We first use the SVD to obtain
its approximated low-rank image with rank $r=50$. The original image and the corresponding approximated low-rank image are displayed in Figure \ref{figure2}. We take $\mathrm{SR}=0.40$ and
$\mathrm{SR}=0.30$ for the low rank image. Two sampled low-rank images with $\mathrm{SR}=0.40$ and $\mathrm{SR}=0.30$ are shown in Figure \ref{figure3}. Numerical results of the three
algorithms for image inpainting are reported in Table \ref{table9}. For $\mathrm{SR}=0.40$ and $\mathrm{SR}=0.30$, we also display recovered Lena image via the three algorithms in
Figure \ref{figure4} and Figure \ref{figure5}, respectively. Comparing these numerical results, we can find that the IGSVT algorithm performs much better than SVT algorithm and ISVT
algorithm on image inpainting for $p=0.5$.

\begin{table}[h!]\footnotesize
\centering
\setlength{\tabcolsep}{0.7mm}{
\begin{tabular}{|c||l|l|l|l|l|l|l|}\hline
\multicolumn{7}{|c|}{SR=0.40}\\\hline
Image&\multicolumn{2}{c}{IGSVT, $p=0.5$}&\multicolumn{2}{|c}{SVT}&\multicolumn{2}{|c|}{ISVT}\\
\hline
(Name,\,rank,\,FR)&RE&Time&RE&Time&RE&Time\\
\hline
(Lena, 50, 2.1531)& 1.38e-05& 43.23& 3.26e-02& 32.93& 1.46e-05& 55.55\\
\hline
\multicolumn{7}{|c|}{SR=0.30}\\\hline
Image&\multicolumn{2}{c}{IGSVT, $p=0.5$}&\multicolumn{2}{|c}{SVT}&\multicolumn{2}{|c|}{ISVT}\\
\hline
(Name,\,rank,\,FR)&RE&Time&RE&Time&RE&Time\\
\hline
(Lena, 50, 1.6149)& 3.02e-05& 159.82& 7.91e-02& 20.96& 3.95e-05& 255.00\\
\hline
\end{tabular}}
\caption{\scriptsize Comparisons of IGSVT algorithm, SVT algorithm and ISVT algorithm for image inpainting}\label{table9}
\end{table}

\section{Conclusions}\label{section6}
It is well known that the affine matrix rank minimization problem is NP-hard and all known algorithms for exactly solving it are doubly exponential in theory and in
practice due to the combinational nature of the rank function. In this paper, inspired by the good performances of the generalized thresholding operator in compressed
sensing, a generalized singular value thresholding operator is generated to solve this NP-hard problem. Numerical experiments on random low-rank matrix completion problems
show that our algorithm performs effectively in finding a low-rank matrix. Moreover, extensive numerical results have illustrated that our algorithm are able to address
low-rank matrix completion problems such as image inpainting. Compared with some state-of-art methods, we can find that our algorithm performs the best on image inpainting.

\section*{Acknowledgments}
The authors would like to thank the reviewers and editors for their useful comments which significantly improve this paper. This research was supported by the
National Natural Science Foundation of China (11771347, 91730306, 41390454, 11271297) and the Science Foundations of Shaanxi Province of China (2016JQ1029, 2015JM1012).





\end{document}